\documentclass[journal]{IEEEtran}

\usepackage{amssymb,amsmath,graphicx,color,psfrag}
\usepackage{enumerate}
\usepackage{url,adjustbox}
\usepackage{graphicx}
\usepackage[utf8]{inputenc}
\usepackage{color,soul}
\usepackage[T1]{fontenc}
\usepackage{upgreek}
\usepackage{multicol}
\usepackage{multirow}
\usepackage{textcomp}
\usepackage{graphicx}
\usepackage{bm}
\usepackage{epstopdf}
\usepackage{amssymb}
\usepackage{amsmath}
\usepackage{psfrag}
\usepackage{amsthm}
\usepackage[tight,footnotesize]{subfigure}
\usepackage{cite}
\usepackage{enumerate}
\usepackage{algorithm}
\usepackage{algpseudocode}
\usepackage{mathtools}
\usepackage{array}
\usepackage{mathabx}
\usepackage{booktabs}
\usepackage{amsmath}
\usepackage{mathrsfs}
\usepackage{array}
\usepackage[tight,footnotesize]{subfigure}
\usepackage{algorithm}
\usepackage{algpseudocode}
\usepackage{multirow}
\usepackage{bm}
\usepackage{mathtools}
\usepackage{soul}

\newcounter{parcount}
\newcounter{parcountA}

\theoremstyle{plain}

\newtheorem{theorem}{Theorem}
\theoremstyle{definition}

\newtheorem{lemma}{Lemma}

\newtheorem{remark}{Remark}

\makeatletter
\newcommand{\vast}{\bBigg@{3.2}}
\newcommand{\Vast}{\bBigg@{5}}

\makeatother

\graphicspath{{Figure/}}

\begin{document}

\title{On Robust Stability and Performance with a Fixed-Order Controller Design for Uncertain Systems }

\author{Jun Ma, 
	Haiyue Zhu, 
	Masayoshi Tomizuka, \IEEEmembership{Life Fellow,~IEEE,}
	and
	Tong Heng Lee 
\thanks{J. Ma and M. Tomizuka are with the Department of Mechanical Engineering, University of California, Berkeley, CA 94720 USA (e-mail:jun.ma@berkeley.edu; tomizuka@berkeley.edu).}
\thanks{H. Zhu is with the Singapore Institute of Manufacturing Technology, A*STAR, Singapore 138634 (e-mail: zhu\_haiyue@simtech.a-star.edu.sg).}
\thanks{T. H. Lee is with the Department of Electrical and Computer Engineering, National University of Singapore, Singapore 117583 (e-mail:eleleeth@nus.edu.sg).}
\thanks{This work has been submitted to the IEEE for possible publication. Copyright may be transferred without notice,	after which this version may no longer be accessible.}
}
\markboth{}
{J. Ma \MakeLowercase{\textit{et al.}}}

\maketitle

\begin{abstract}
Typically, it is desirable to design a control system that is not only robustly stable in the presence of parametric uncertainties but also guarantees an adequate level of system performance. However, most of the existing methods need to take all extreme models over
an uncertain domain into consideration,
which then results in costly computation.
Also, since these approaches
attempt rather unrealistically to
guarantee the system performance over a full frequency range, a conservative design is always admitted.
Here, taking a specific viewpoint of robust stability and performance under a stated restricted frequency range
(which is applicable in rather many real-world situations),
this paper provides an essential basis for the design of a fixed-order controller for a system with bounded parametric uncertainties, which avoids the tedious but necessary evaluations
of the specifications on all the extreme models. A Hurwitz polynomial is used in the design and the robust stability is characterized by the notion of positive realness,
such that the required robust stability condition is then suitably successfully constructed.
Also, the robust performance criteria in terms of sensitivity shaping under different frequency ranges are constructed based on
an approach of bounded realness analysis. Furthermore, the conditions for robust stability and performance are expressed in the framework of linear matrix inequality (LMI) constraints,
and thus
can be efficiently solved. Comparative simulations are provided to demonstrate the effectiveness and efficiency of the proposed approach.
\end{abstract}


\section{Introduction}
In quite a large number of control engineering applications,
the presence of parametric uncertainties have been challenging issues to deal with,
and there is a clear and evident need to consider
their effects on the deterioration of system performance or even the cause of instability.
There thus have been many great efforts
which are devoted to
developments in
robust control for designing a control system that guarantees the
required robust stability~\cite{celentano2018new}.
Among all the techniques for robust stabilization,
quadratic stabilization theory
provides an effective basis to cater to
these parametric uncertainties~\cite{petersen1987stabilization,zhou1988robust}.
Related work in \cite{petersen1987disturbance} presents a controller design approach to guarantee a prespecified disturbance attenuation level based on the algebraic Riccati equation.
Further, these results have been extended in~\cite{khargonekar1988h} to solve the $H_\infty$ control problem.
As a methodology with attractive applicability,
$H_\infty$ control has been extensively researched and developed over the
period of the 1980s,
particularly in view of the rather commonly needed requirement
of robust stabilization and disturbance attenuation.
In this very attractive approach,
the frequency-domain characteristics
are also expressed
by related mathematical statements in the
time domain.
With such a framework, developments
proceeded such that
in~\cite{doyle1989state}, several fundamental results on $H_\infty$ control are established.
However, these results are only for the classes of well-known systems.
Then in~\cite{khargonekar1990robust}, a linear system with norm-bounded parametric uncertainties is related to ARE-based $H_\infty$ norm conditions.
Moreover, robust $H_\infty$ performance problems are addressed in~\cite{xie1990robust,xie1992robust,de1993h,gu1994h,xie1996output} and a class of systems with norm-bounded parametric uncertainties
is also taken into consideration.


Substantial further works on guaranteed cost control have also been presented
~\cite{petersen1994optimal,kosmidou2006linear,ma2020arxiv,ma2020symmetric,guo2014robust},
which aim to design a control system such that an upper bound of the quadratic performance is guaranteed for all admissible parametric uncertainties.
Additionally in~\cite{peres1993h}, a stabilizing controller is proposed for a class of systems with convex-bounded parametric uncertainties, which assures a certain specified attenuation level.
However, to obtain a deterministic and reliable solution,
it is required to check the stability and performance of the vertex set consisting of several extreme matrices~\cite{geromel1991convex}.
Because the number of extreme systems to be checked increases exponentially with the uncertain parameters,
the computation is always costly~\cite{song2017robust,ma2019robust,ma2019parameter}.
Meanwhile, other extensive developments have been made
on the robustness property attainment involving interval matrix uncertainties
to reduce the number of vertices that are required to be checked.
For example, the work in \cite{alamo2008new} presents a new vertex result for the robust synthesis problem.
But here too, if the interval uncertainties appear in the linear matrix inequality (LMI) in an affine way,
a rather substantial number of LMIs will need to be integrated to be solved, and the solution is still exponential.
Yet in some other scenarios, a randomized algorithm approach is employed to ensure
the required stochastic robust stability and performance.
However, most of these probabilistic approaches are solved based on Monte-Carlo simulation~\cite{guo2015efficient},
and they are typically not considered as practically preferred due to these excessive simulations.

In view of the great potential for practical usage,
the methodology of fixed-order controllers
has attracted considerable attention for finite frequency specifications due to their simplicity, reliability,
and ease of implementation~\cite{karimi2018robust}.
Though the design of a fixed-order controller in the presence of parametric uncertainties is NP-hard,
a variety of design approaches have been rather successfully used
including bilinear matrix inequality (BMI)~\cite{safonov1994control,goh1994global},
convex approximation~\cite{karimi2007robust,khatibi2008fixed} and iterative heuristic optimization~\cite{iwasaki1999dual}.
It is worthwhile to mention that for all these above-mentioned control methods,
the robust performance property is achieved over the full frequency range.
However, in many real-world situations,
the control performance specifications
are typically
only specified
and of pertinent interest
within a stated frequency range
for the real-world system to be controlled~\cite{li2013heuristic,li2019enhanced}.
Therefore, the design from the perspective of the full frequency range is overly conservative.
To cater to the more realistic practical requirement of a restricted frequency range,
recent researches reveal that these control strategies
can be combined with certain appropriate frequency weighting functions~\cite{bevrani2015robust};
but here, a required concomitant
strong computational capability of the hardware is a burden in practical usage
because these weighting functions invariably cause a marked increase in the system orders.
Under these circumstances, a rather significant and important open problem therefore exists
on an effective and efficient approach to design a fixed-order controller
for a system affected by bounded parametric uncertainties,
and which considers robust stability and performance under the condition of a stated restricted frequency range,
but without all the afore-mentioned practical difficulties.


Thus in this work, the key objective is to design a fixed-order controller for an uncertain system
that guarantees closed-loop stability and a suitably adequate level of performance under a stated restricted frequency range,
and which addresses all the short-comings previously highlighted.
As part of the methodology in the development in this work,
the connection between the time-domain method and frequency-domain specifications are established
to meet the stability criteria and the performance specifications.
Several theorems are then developed and provided
to support the derivation of the main results.
Then, an LMI condition is given
to ensure the robust stability of the closed-loop system in the presence of the bounded parametric uncertainties.
Additionally further, certain LMI conditions are also provided
for the robust performance of the system which is guaranteed in terms of sensitivity shaping. Consequently, the required controller is obtained by solving the above-mentioned LMIs.


The remainder of this paper is organized as follows.
In Section II, the necessary preliminaries on the closed-loop control of a system with bounded parametric uncertainties and frequency range characterization are provided.
Then, a first set of newly developed theoretical results in this work on the robust stability condition is presented in Section III.
Next, further theoretical results developed in this work on robust performance criteria are presented in Section IV.
In Section V, to validate the new proposed controller design approach,
a numerical example is provided with simulation results to show its effectiveness.
Finally, conclusions are drawn in Section VI.


\textbf{Notations:} For matrix $A$, the symbols $A^T$ and $A^*$ represent the transpose and the complex conjugate transpose of a matrix, respectively. $\text{Re}(A)$ denotes the real part of a matrix. $I$ represents the identity matrix with appropriate dimensions. $\textup{diag}\{a_1, a_2, \cdots, a_n\}$ represents the diagonal matrix with numbers $a_1, a_2,  \cdots, a_n$ as diagonal entries. $\mathbb{R}$ and $\mathbb{C}$ indicate the sets of real and complex matrices, respectively. $\mathbb{H}_n$ stands for the set of $n \times n$ complex Hermitian matrices. $f \ast g$ denotes the convolution operation of two functions. The operator $\otimes$ represents the Kronecker's product. 
The symbol $s$ in the bracket, when it typically appears in such expressions as $T(s)$, etc., represents the Laplace variable.

\section{Preliminaries}
\subsection{Problem statement}~\label{sec:2.1}
As in typical nomenclature, the Single-Input-Single-Output (SISO) plant is represented by an $n$th-order rational transfer function in continuous time, and is given by
\begin{equation}\label{PlantModel1}
\begin{aligned}
P(s)=\frac{b_{1}s^{n-1}+\cdots+b_{n}}{s^{n}+a_{1}s^{n-1}+\cdots+a_{n}},
\end{aligned}
\end{equation}
where $a_{i}$ and $b_{i}$ are uncertain parameters with $a_{i}\in[a_{i}^{l},a_{i}^{u}]$ and $b_{i}\in[b_{i}^{l},b_{i}^{u}]$, $i=1,2,\cdots,n$.


Next, define $a_{i}^{c}=(a_{i}^{l}+a_{i}^{u})/2$ and $b_{i}^{c}=(b_{i}^{l}+b_{i}^{u})/2$ as the medians of the uncertain parameters $a_{i}$ and $b_{i}$, respectively. Similarly, further define $a_{i}^{d}=(a_{i}^{u}-a_{i}^{l})/2$ and $b_{i}^{d}=(b_{i}^{u}-b_{i}^{l})/2$ as the deviations of the uncertain parameters $a_{i}$ and $b_{i}$, respectively. Then, $a_{i} = a_{i}^{c}+ a_{i}^d \delta_{ai}$, $b_{i} = b_{i}^{c}+ b_{i}^d \delta_{bi}$, where $\delta_{ai}\in [-1,1]$ and $\delta_{bi}\in [-1,1]$ are standard interval variables.

Also, define $\Delta_{a}=\textup{diag}\{\delta_{a1}, \delta_{a2}, \cdots, \delta_{an}\}$ and $\Delta_{b}=\textup{diag}\{\delta_{b1}, \delta_{b2}, \cdots, \delta_{bn}\}$, and the plant \eqref{PlantModel1} can be expressed by
\begin{equation}\label{PlantModel2}
\begin{aligned}
P(s)=\frac{(\bm{b^{c}}+[0\quad \bm{b_{d}}\Delta_{b}])\bm{s_{n}}^{T}}{(\bm{a^{c}}+[0\quad \bm{a_{d}}\Delta_{a}])\bm{s_{n}}^{T}},
\end{aligned}
\end{equation}
where
$\bm{a^{c}}=[1\quad a_{1}^{c}\quad a_{2}^{c}\quad \cdots\quad a_{n}^{c}]$, $\bm{b^{c}}=[0\quad b_{1}^{c}\quad b_{2}^{c}\quad \cdots\quad b_{n}^{c}]$, $\bm{a_{d}}=[a_{1}^{d}\quad a_{2}^{d}\quad \cdots\quad a_{n}^{d}]$, $\bm{b_{d}}=[b_{1}^{d}\quad b_{2}^{d}\quad \cdots\quad b_{n}^{d}]$, $\bm{s_{n}}=[s^{n}\quad s^{n-1}\quad \cdots\ s\quad 1]$. 	For brevity, define $\bm{a}=\bm{a^{c}}+[0\quad \bm{a_{d}}\Delta_{a}]$ and $\bm{b}=\bm{b^{c}}+[0\quad \bm{b_{d}}\Delta_{b}]$.


By the standard negative feedback configuration, an $m$th-order controller is to be designed, which is given by
\begin{equation}\label{Controller1}
\begin{aligned}
K(s)=\frac{y_{0}s^{m}+y_{1}s^{m-1}+\cdots+y_{m}}{s^{m}+x_{1}s^{m-1}+\cdots+x_{m}}.
\end{aligned}
\end{equation}
Equivalently, \eqref{Controller1} is expressed by
\begin{equation}\label{Controller2}
\begin{aligned}
K(s)=\frac{\bm{y}\bm{s_{m}}^{T}}{\bm{x}\bm{s_{m}}^{T}},
\end{aligned}
\end{equation}
where
$\bm{x}=[1\quad x_{1}\quad x_{2}\quad \cdots\quad x_{m}]$,
$\bm{y}=[y_{0}\quad y_{1}\quad y_{2}\quad \cdots \quad y_{m}]$,
$\bm{s_{m}}=[s^{m}\quad s^{m-1}\quad \cdots \quad s\quad 1]$.


For the closed-loop system, the sensitivity transfer function $S(s)$ and the complementary sensitivity transfer function $T(s)$ are given by
\begin{equation}\label{SFunction}
S(s)=\frac{S_{num}}{S_{den}}, \quad T(s)=\frac{T_{num}}{T_{den}},
\end{equation}
respectively, where
$S_{num} = (\bm{a}\ast\bm{x})\bm{s_{m+n}}^{T}$,
$T_{num} = (\bm{b}\ast\bm{y})\bm{s_{m+n}}^{T}$,
$S_{den}=T_{den} = (\bm{a}\ast\bm{x}+\bm{b}\ast\bm{y})\bm{s_{m+n}}^{T}$,
with
$
\bm{s_{m+n}}=[s^{m+n}\quad s^{m+n-1}\quad \cdots \quad s\quad 1].
$
Equivalently, we have
\begin{equation}
\begin{aligned}
S_{num} &= (\bm{a^{c}}\ast\bm{x})\bm{s_{m+n}}^{T}+((\bm{a_{d}}\Delta_{a})\ast\bm{x})\bm{s_{m+n-1}}^{T},\\
T_{num} &= (\bm{b^{c}}\ast\bm{y})\bm{s_{m+n}}^{T}+((\bm{b_{d}}\Delta_{b})\ast\bm{y})\bm{s_{m+n-1}}^{T},\\
S_{den}&=T_{den}=(\bm{a^{c}}\ast\bm{x}+\bm{b^{c}}\ast\bm{y})\bm{s_{m+n}}^{T}\\&\qquad\qquad\,\,+((\bm{a_{d}}\Delta_{a})\ast\bm{x}+(\bm{b_{d}}\Delta_{b})\ast\bm{y})\bm{s_{m+n-1}}^{T},
\end{aligned}
\end{equation}
where $\bm{s_{m+n-1}}=[s^{m+n-1}\quad s^{m+n-2}\quad \cdots\quad s\quad 1].$

Thus the objective here is to design a fixed-order controller $K(s)$ for the uncertain system $P(s)$ such that\\
(1). The robust stability of the closed-loop system is guaranteed in the presence of parametric uncertainties.\\
(2). The robust performance specifications of the closed-loop system in terms of sensitivity shaping are satisfied, i.e. $\big| S(j\omega)\big| < \rho_s$, $\omega \in \Omega_s$ and $\big| T(j\omega)\big| < \rho_t$, $\omega \in \Omega_t$.

\begin{remark}
	In the literature, many instances of controller design methods result in the controller with a larger or equal order as that of the plant. However, this might be a restrictive condition in certain scenarios because the implementation of such controllers will lead to high cost and fragility. On the other hand, there has been a considerable interest in design of low-order controllers to facilitate the practical implementation. Normally the selection of controller order is decided by the user in view of specific situations.
\end{remark}

\subsection{Frequency range characterization}
Basically, a frequency range can be visualized as a curve on the complex plane.
Note that a curve on the complex plane is a collection of points $\lambda(t) \in \mathbb{C}$ continously parameterized by $t$, for $t_0\leq t \leq t_f$, where $t_0, t_f \in \mathbb{R}\cup \{\pm\infty\}$,  which can be	characterized by a set~\cite{iwasaki2005generalized} 	
\begin{equation}
\Lambda = \{\lambda \in \mathbb{C}: \sigma(\lambda,\Phi)=0, \sigma(\lambda,\Psi) \geq 0\},
\end{equation}
where $\Phi$, $\Psi\in \mathbb{H}_2$ and
\begin{equation}
\sigma(\lambda,\Phi)=
\left[\begin{array}{c}
\lambda\\1
\end{array}
\right]^* \Phi
\left[\begin{array}{c}
\lambda\\1
\end{array}\right]
,\quad
\sigma(\lambda,\Psi)=
\left[\begin{array}{c}
\lambda\\1
\end{array}\right]^* \Psi
\left[\begin{array}{c}
\lambda\\1
\end{array}\right].
\end{equation}
For the continuous time domain, one has
\begin{equation}
\Phi =
\left[\begin{array}{cc}
0 & 1\\1& 0
\end{array}\right],\quad
\Lambda=\{j\omega: \omega\in\Omega\},
\end{equation}
where $\Omega$ is a subset of real numbers specified by $\Psi$. Table 1 summarizes the characterization of finite frequency ranges, and more details can be referred in~\cite{iwasaki2005generalized} .

\begin{table}[h]~\label{tab:freq}\centering
	\caption{Characterization of frequency range}
	\begin{adjustbox}{width=\columnwidth,center}
	\begin{tabular}{cccc}
		\toprule
		$\Omega$ & $\omega\in[0,\,\omega_l]$ & $\omega\in[\omega_l,\,\omega_h]$ & $\omega\in[\omega_h,\,+\infty)$ \\
		\midrule
		$\Psi$ &
		$\left[
		\begin{array} {cc}
		-1 & 0\\
		0 &  \omega_l^2
		\end{array}
		\right]$
		&
		$\left[
		\begin{array} {cc}
		-1 & j\frac{\omega_h+\omega_l}{2}\\
		-j \frac{\omega_h+\omega_l}{2} & -\omega_l\omega_h
		\end{array}
		\right]$
		&
		$\left[
		\begin{array} {cc}
		1 & 0\\
		0 & -\omega_h^2
		\end{array}
		\right]$
		\\
		\bottomrule
			\end{tabular}
		\end{adjustbox}
\end{table}

\subsection{Lemmas used in the sequel}

In this subsection, the following lemmas are presented to be used in the sequel for the derivation of the main results.

\begin{lemma}~\cite{xie1996output}\label{lemma:xie}
	Given matrices $Q$, $H$, $E$ and $R$ of appropriate dimensions and with $Q$ and $R$ symmetrical and $R>0$, then
	\begin{equation}
	Q+HFE+E^TF^TH^T<0,
	\end{equation}
	for all $F$ satisfying $F^TF\leq R$, if and only if there exists some $\varepsilon >0$ such that
	\begin{equation}
	Q+\varepsilon^2 HH^T+\varepsilon^{-2}E^T RE<0.
	\end{equation}
\end{lemma}

\begin{lemma}\cite{sun1994solution}~\label{lemma:pos}
	Consider $\sum \triangleq \{A,B,C,D\}$ as a minimal state-space realization of a rational and proper transfer function $G(s)$, the positive realness condition
	\begin{equation}\label{}
	\begin{aligned}
	\textup{Re}\big(G(s)\big)>0
	\end{aligned},
	\end{equation}
	is guaranteed if and only if there admits a Hermitian matrix $P>0$ such that\\
	\begin{equation}\label{KYPPR}
	\begin{aligned}
	\left[
	\begin{array} {cc}
	A & B  \\
	I & 0
	\end{array}
	\right]^{T}
	\left[
	\begin{array} {cc}
	0 & P  \\
	P & 0
	\end{array}
	\right]
	\left[
	\begin{array} {cc}
	A & B  \\
	I & 0
	\end{array}
	\right]-
	\left[
	\begin{array} {cc}
	0 & C^{T}  \\
	C & D+D^{T}
	\end{array}
	\right]<0
	\end{aligned}.
	\end{equation}
\end{lemma}

\begin{lemma}~\cite{iwasaki2005generalized}~\label{lemma: GKYP}
	Consider $(A, B, C, D)$ as a minimal state-space realization of a rational transfer function $G(s)$, given $\rho>0$, the finite frequency bounded realness condition
	\begin{equation}\label{}
	\begin{aligned}
	\big| G(j\omega)\big|<\rho,\ \omega\in\Omega
	\end{aligned},
	\end{equation}
	is guaranteed if and only if there exist Hermitian matrices $P$ and $Q>0$, such that
	\begin{equation}\label{GKYPbg}
	\begin{aligned}
	\left[
	\begin{array} {cc}
	\left[
	\begin{array} {cc}
	A & B  \\
	I & 0
	\end{array}
	\right]^{T}
	\Upxi
	\left[
	\begin{array} {cc}
	A & B  \\
	I & 0
	\end{array}
	\right]
	+
	\left[
	\begin{array} {cc}
	0 & 0  \\
	0 & -\rho
	\end{array}
	\right]
	&
	\begin{array} {c}
	C^{T} \\
	D^{T}
	\end{array}
	\\
	\begin{array} {cc}
	C & D
	\end{array}
	&
	-\rho
	\end{array}
	\right]<0
	\end{aligned},
	\end{equation}
	where $\Upxi=\Phi\otimes P+\Psi\otimes Q$, and $\Phi$ and $\Psi$ are matrices to characterize the frequency range $\Omega$.
\end{lemma}

\section{Robust stability characterization via real positiveness analysis}

In this section, first of all, pertinent new results Theorem~\ref{theorem:1} and Theorem~\ref{theorem:2} are presented, which are used in the sequel for deriving the robust stability condition.
\begin{theorem}~\label{theorem:1}
	Given matrices $Q$, $H$, $E$ with appropriate dimensions, $Q$ is symmetrical, $\Delta=\textup{diag}\{\delta_1, \delta_2, \cdots, \delta_n\}$ with $\delta_i\in [-1,1]$, $i=1,2,\cdots,n$,
	\begin{equation}
	Q+\left[\begin{array}{cc}
	0 & E^T\Delta H^T\\
	H\Delta E &0
	\end{array}\right]<0 \label{eqn:theorem1:1}
	\end{equation}
	holds if and only if there exists a matrix $R=\textup{diag} \{\varepsilon_1 , \varepsilon_2, \cdots, \varepsilon_n\}$ with $\varepsilon_i>0$, $i=1,2,\cdots,n$, such that
	\begin{equation}
	Q+
	\left[
	\begin{array}{cc}
	E^T R^{-1}E & 0\\
	0& H R H^T
	\end{array}
	\right] <0. \label{eqn:theorem1:2}
	\end{equation}
\end{theorem}
\noindent{\textbf{Proof of Theorem~\ref{theorem:1}:}} Sufficiency: By matrix decomposition, \eqref{eqn:theorem1:1} is equivalent to
\begin{equation}~\label{eq:decom1}
Q+\left[\begin{array}{c}0\\H\end{array}\right]\Delta\left[\begin{array}{cc}E&0\end{array}\right]+\left[\begin{array}{cc}E&0\end{array}\right]^T\Delta\left[\begin{array}{c}
0\\H \end{array}\right]^T <0.
\end{equation}	
$\begin{bmatrix} 0&H\end{bmatrix}^T$ can be partitioned as $n$ column vectors, i.e. $H_{1}, H_{2}, \cdots, H_{n}$; similarly, $\begin{bmatrix}E&0\end{bmatrix}$ can be partitioned as $n$ row vectors, i.e. $E_{1}, E_{2}, \cdots, E_{n}$. Then, it is easy to obtain
\begin{equation}
\begin{aligned}
&\quad \left[\begin{array}{c}0\\H\end{array}\right]\Delta \left[\begin{array}{cc}
E&0 \end{array}\right]+\left[\begin{array}{cc}E&0\end{array}\right]^T\Delta\left[\begin{array}{c}
0\\H \end{array}\right]^T \nonumber\\&= \sum_{i=1}^n(H_i\delta_iE_i+E_i^T\delta_iH_i^T).
\end{aligned}
\end{equation}
Thus, \eqref{eq:decom1} is equivalent to
\begin{equation}~\label{eq:decom2}
Q+\sum_{i=1}^n(H_i\delta_iE_i+E_i^T\delta_iH_i^T)<0.
\end{equation}
Define
\begin{equation}~\label{eqn:Qn-1}
Q_{n-1}=Q+\sum \limits_{i=1}^{n-1}(H_i\delta_iE_i+E_i^T\delta_iH_i^T),
\end{equation} then \eqref{eq:decom2} can be written as
\begin{equation}~\label{eq:decom3}
Q_{n-1}+H_n\delta_nE_n+E_n^T\delta_nH_n^T<0.
\end{equation}
From Lemma~\ref{lemma:xie}, \eqref{eq:decom3} holds if and only if there exists $\varepsilon_n>0$ such that
\begin{equation}~\label{eqn:18}
Q_{n-1}+\varepsilon_nH_n H_n^T+\varepsilon_n^{-1} E_n^T E_n<0.
\end{equation}
Similarly, define
\begin{equation}~\label{eqn:Qn-2}
Q_{n-2}=Q+\sum \limits_{i=1}^{n-2}(H_i\delta_iE_i+E_i^T\delta_iH_i^T)+\varepsilon_n H_n H_n^T+\varepsilon_n^{-1}E_n^T E_n,\end{equation}
then we have
\begin{equation}~\label{eq:decom4}
Q_{n-2}+H_{n-1}\delta_{n-1}E_{n-1}+E_{n-1}^T\delta_{n-1}H_{n-1}^T<0.
\end{equation}
Again, from Lemma~\ref{lemma:xie}, \eqref{eq:decom4} holds if and only if there exists $\varepsilon_{n-1}>0$ such that
\begin{equation}\label{eqn:20}
Q_{n-2}+\varepsilon_{n-1}H_{n-1} H_{n-1}^T+\varepsilon_{n-1}^{-1} E_{n-1}^T E_{n-1}<0.
\end{equation}
Substitute~\eqref{eqn:Qn-2} to~\eqref{eqn:20}, we have
\begin{equation}
\begin{aligned}
&Q+\sum \limits_{i=1}^{n-2}(H_i \delta_i E_i+E_i^T \delta_i^{-1}H_i^T)+\varepsilon_{n-1}H_{n-1}H_{n-1}^T \\
&+\varepsilon_{n-1}^{-1}E_{n-1}^TE_{n-1} + \varepsilon_{n}H_{n}H_{n}^T+\varepsilon_n^{-1}E_{n}^TE_{n} <0.
\end{aligned}
\end{equation}

In a similar way, it is straightforward that \eqref{eq:decom1} holds if and only if there exist $\varepsilon_1, \varepsilon_2, \cdots, \varepsilon_n>0$ such that
\begin{equation}
Q+\sum \limits_{i=1}^n (\varepsilon_i H_i H_i^T+\varepsilon_i^{-1}  E_i^T E_i)<0.	
\end{equation}
Define  $R=\textup{diag}\{\varepsilon_1, \varepsilon_2, \cdots, \varepsilon_n\}$, we have
\begin{equation}
Q+\left[\begin{array}{cc}E&0\end{array}\right]^TR^{-1}\left[\begin{array}{cc}E&0\end{array}\right]+\left[\begin{array}{c}0\\H\end{array}\right] R\left[\begin{array}{c}
0\\H \end{array}\right]^T <0,
\end{equation}
which can be further expressed by
\begin{equation}
Q+
\left[
\begin{array}{cc}
E^T R^{-1}E & 0\\
0& H R H^T
\end{array}
\right] <0.
\end{equation}
\\
Necessity: it can be proved in a similar way thus is omitted.
\hfill{\qed}

\begin{theorem}~\label{theorem:2}
	Given matrices $Q$, $H_i$, $E_i$, $i=1,2,\cdots,m$   with appropriate dimensions, $Q$ is symmetrical, $\Delta_i=\textup{diag}\{\delta_{i1}, \delta_{i2}, \cdots, \delta_{in}\}$ with $\delta_{ij}\in [-1,1]$, $i=1,2,\cdots, m$, $j=1,2,\cdots,n$,
	\begin{equation}
	Q+
	\left[\begin{array}{cccc}
	0 & 	\sum \limits_{i=1}^m E_i^T\Delta_i H_i^T\\
	\sum \limits_{i=1}^m H_i\Delta_i E_i &0
	\end{array}\right]<0 \label{eqn:theorem3:1}
	\end{equation}
	holds if and only if there exist matrices $R_i=\textup{diag}\{\varepsilon_{i1}, \varepsilon_{i2}, \cdots, \varepsilon_{in}\}$ with $\varepsilon_{ij}>0$, $i=1,2,\cdots, m$, $j=1,2,\cdots,n$, such that
	\begin{equation}
	Q+
	\left[
	\begin{array}{cc}
	\sum \limits_{i=1}^m E_i^T R_i^{-1}	E_i & 0\\
	0& \sum \limits_{i=1}^m H_i  R_i H_i^T \end{array}
	\right]<0. \label{eqn:theorem2:2}
	\end{equation}
\end{theorem}
\noindent{\textbf{Proof of Theorem~\ref{theorem:2}:}}
Theorem~\ref{theorem:2} is essentially an extension of Theorem~\ref{theorem:1}, the proof is rather straightforward and is omitted here. \hfill{\qed}

Next, a Hurwitz polynomial $d_c(s)$ is selected and an associated transfer function is defined as
\begin{equation}\label{DefGs}
\begin{aligned}
G_{s}(s)=&\frac{(\bm{a}\ast\bm{x}+\bm{b}\ast\bm{y})\bm{s_{m+n}}^{T}}{d_{c}(s)}\\
=&G_{sn}(s)+G_{su}(s),
\end{aligned}
\end{equation}
where $G_{sn}(s)$ and $G_{su}(s)$ are the nominal part and uncertain part of $G_{s}(s)$, respectively, which are given by
\begin{equation}\label{DefGsnp}
\begin{aligned}
G_{sn}(s)=&\frac{(\bm{a^{c}}\ast\bm{x}+\bm{b^{c}}\ast\bm{y})\bm{s_{m+n}}^{T}}{d_{c}(s)},\\
G_{su}(s)=&\frac{((\bm{a_{d}}\Delta_{a})\ast\bm{x}+(\bm{b_{d}}\Delta_{b})\ast\bm{y})\bm{s_{m+n-1}}^{T}}{d_{c}(s)}.
\end{aligned}
\end{equation}

The realization of  $G_{sn}(s)$ in the controllable canonical form is denoted by
\begin{equation}\label{realization1}
\sum\nolimits_{sn} \triangleq
\{A_{sn}, B_{sn},C_{sn}, D_{sn}\},
\end{equation}
and then the realization of $G_s(s)$ is given by
\begin{equation}~\label{eq:realization}
\sum\nolimits_{s}  \triangleq
\{
A_{sn}, B_{sn},
C_{sn}+\bm{a_{d}}\Delta_{a}X+\bm{b_{d}}\Delta_{b}Y, D_{sn}
\},
\end{equation}
where 		$X$ and $Y$ are the Toeplitz matrices given by
\begin{equation}\label{RealNominal}
\begin{aligned}
X=&
\left[
\begin{array} {ccccccc}
1       & x_{1}   & \cdots    & x_{m}     &0        &0          &0      \\
0       & 1       & x_{1}     & \cdots    & x_{m}   &0          &0      \\
\vdots  &\vdots   &\vdots     &\ddots     &\vdots   &\vdots     &\vdots\\
0       &0        &0          & 1       & x_{1}     & \cdots    & x_{m}
\end{array}
\right],\\
Y=&
\left[
\begin{array} {ccccccc}
y_{0}   & y_{1}   & \cdots    & y_{m}     &0        &0          &0      \\
0       & y_{0}   & y_{1}     & \cdots    & y_{m}   &0          &0      \\
\vdots  &\vdots   &\vdots     &\ddots     &\vdots   &\vdots     &\vdots\\
0       &0        &0          & y_{0}     & y_{1}   & \cdots    & y_{m}
\end{array}
\right].\\
\end{aligned}
\end{equation}
Notably here, the stability of the system is closely related to the notion of positive realness, and it is guaranteed if and only if
\begin{equation}\label{StabilityPR}
\begin{aligned}
\text{Re}\big(G_{s}(s)\big)>0.
\end{aligned}
\end{equation}
Then, to construct a condition for the needed robust stability in the presence of parametric uncertainties, Theorem~\ref{theorem:stability} is proposed.

\begin{theorem}~\label{theorem:stability}
	The robust stability of the system~\eqref{PlantModel1} in the presence of bounded parametric uncertainties characterized by standard interval variables is guaranteed under the controller~\eqref{Controller1} if and only if there exist a Hermitian matrix $P_s>0$, diagonal matrices $R_{sa}>0$ and $R_{sb}>0$ such that
	\begin{equation}\label{LMI1}
	\begin{aligned}
	\left[
	\begin{array} {c|c}
	\Gamma_s  &
	\begin{array} {cc}
	X^T & Y^T \\
	0 & 0
	\end{array}  \\ \hline
	\begin{array} {cc}
	X & 0 \\
	Y & 0
	\end{array}
	&
	\begin{array} {cc}
	-R_{sa} & 0  \\
	0 &  -R_{sb}
	\end{array}
	\end{array}
	\right]
	<0
	\end{aligned},
	\end{equation}
	where
	\begin{equation}
	\begin{aligned}
	\Gamma_s &=
	\left[
	\begin{array} {cc}
	A_{sn} & B_{sn}  \\
	I & 0
	\end{array}
	\right]^{T}
	\left[
	\begin{array} {cc}
	0 & P_s  \\
	P_s & 0
	\end{array}
	\right]
	\left[
	\begin{array} {cc}
	A_{sn} & B_{sn}  \\
	I & 0
	\end{array}
	\right]\\
	&-\left[
	\begin{array} {cc}
	0 & C_{sn}^T  \\
	C_{sn} & D_{sn}+D_{sn}^T-\bm{a_{d}} R_{sa} \bm{a_{d}}^T-\bm{b_{d}} R_{sb} \bm{b_{d}}^T
	\end{array}
	\right].
	\end{aligned}
	\end{equation}
\end{theorem}

\noindent{\textbf{Proof of Theorem~\ref{theorem:stability}:}} From Lemma~\ref{lemma:pos}, considering the state-space realization~\eqref{eq:realization},  the positive realness condition~\eqref{StabilityPR} is guaranteed if and only if there exists a Hermitian matrix $P_{s}>0$ such that\\
\begin{equation}\label{KYPPR1}
\begin{aligned}
&\left[
\begin{array} {cc}
A_{sn} & B_{sn}  \\
I & 0
\end{array}
\right]^{T}
\left[
\begin{array} {cc}
0 & P_s  \\
P_s & 0
\end{array}
\right]
\left[
\begin{array} {cc}
A_{sn} & B_{sn}  \\
I & 0
\end{array}
\right]\\
&-\left[
\begin{array} {cc}
0 & C_{sn}^T \\
C_{sn}& D_{sn}+D_{sn}^{T}
\end{array}
\right] \\
&+\left[
\begin{array} {cc}
0 &  -X^T  \Delta_{a}  \bm{a_{d}}^T  -Y^T \Delta_{b}  \bm{b_{d}}^T \\
-\bm{a_{d}} \Delta_{a} X-\bm{b_{d}} \Delta_{b} Y & 0
\end{array}
\right]<0.
\end{aligned}
\end{equation}


From Theorem~\ref{theorem:2}, \eqref{KYPPR1} holds if and only if there exist positive definite diagonal matrices $R_{sa}$ and $R_{sb}$ such that
\begin{equation}
\begin{aligned}
&\left[
\begin{array} {cc}
A_{sn} & B_{sn}  \\
I & 0
\end{array}
\right]^{T}
\left[
\begin{array} {cc}
0 & P_s  \\
P_s & 0
\end{array}
\right]
\left[
\begin{array} {cc}
A_{sn} & B_{sn}  \\
I & 0
\end{array}
\right]\\
&-
\left[
\begin{array} {cc}
0 & C_{sn}^T  \\
C_{sn} & D_{sn}+D_{sn}^T
\end{array}
\right]\\
& +\left[
\begin{array} {cc}
X^T R_{sa}^{-1} X +Y^T R_{sb}^{-1} Y   & 0 \\
0&  \bm{a_{d}} R_{sa} \bm{a_{d}}^{T}+ \bm{b_{d}}R_{sb}\bm{b_{d}}^{T}
\end{array}
\right]
<0.
\end{aligned}
\end{equation}
Then, we have
\begin{equation}\label{KYPPR2}
\begin{aligned}
&\left[
\begin{array} {cc}
A_{sn} & B_{sn}  \\
I & 0
\end{array}
\right]^{T}
\left[
\begin{array} {cc}
0 & P_s  \\
P_s & 0
\end{array}
\right]
\left[
\begin{array} {cc}
A_{sn} & B_{sn}  \\
I & 0
\end{array}
\right]\\
&-
\left[
\begin{array} {cc}
0 & C_{sn}^T  \\
C_{sn} & D_{sn}+D_{sn}^T
\end{array}
\right]\\
& +\left[
\begin{array} {c}
0 \\
\bm{a_{d}}
\end{array}
\right]R_{sa}
\left[
\begin{array} {cc}
0 &  \bm{a_{d}}^{T}
\end{array}
\right]
+
\left[
\begin{array} {c}
0 \\
\bm{b_{d}}
\end{array}
\right]R_{sb}
\left[
\begin{array} {cc}
0 &  \bm{b_{d}}^{T}
\end{array}
\right]
\\
&+\left[
\begin{array} {c}
X^T \\
0
\end{array}
\right] R_{sa}^{-1}
\left[
\begin{array} {cc}
X & 0
\end{array}
\right]
+
\left[
\begin{array} {c}
Y^T \\
0
\end{array}
\right] R_{sb}^{-1}
\left[
\begin{array} {cc}
Y & 0
\end{array}
\right]
<0,
\end{aligned}
\end{equation}
which can be further expressed in the form of \eqref{LMI1}. This completes the proof of the theorem. \hfill{\qed}



\section{Robust performance characterization with sensitivity shaping}
To develop an effective methodology for robust performance characterization with sensitivity shaping,
it is pertinent to note here that
robust performance specifications are suitably characterized,
where the infinity norm of the sensitivity function and the complementary sensitivity function are bounded by certain values. With a given Hurwitz polynomial,
the stated bound condition on the infinity norm of a rational transfer function can be separated by two conditions. To summarize this finding, Theorem~\ref{theorem:convex} is now presented.
\begin{theorem}~\label{theorem:convex}
	Consider a rational transfer function $G(s)=n(\lambda,s)/d(\lambda,s)$, where $\lambda$ is a parameter vector appeared affinely in the polynomials $n(\lambda,s)$ and $d(\lambda,s)$. For any given Hurwitz polynomial $d_{c}(s)$, $\big|G(s)\big|<\rho$ is guaranteed if the following two conditions hold:\\
	\begin{equation}\label{eq:convex_1}
	\begin{aligned}
	\left|\frac{n(\lambda,s)}{d_{c}(s)}\right|<(1-\delta)\rho
	\end{aligned},
	\end{equation}
	and
	\begin{equation}\label{eq:convex_2}
	\begin{aligned}
	\left|1-\frac{d(\lambda,s)}{d_{c}(s)}\right|<\delta
	\end{aligned},
	\end{equation}
	where $\delta\in(0,1)$.
\end{theorem}

\noindent{\textbf{Proof of Theorem~\ref{theorem:convex}:}}
From \eqref{eq:convex_2}, we have
\begin{equation}\label{eq:convex_3}
1-\delta< \left|  \frac{d(\lambda,s)}{d_{c}(s)}\right| < 1+\delta.
\end{equation}
By \eqref{eq:convex_1}/\eqref{eq:convex_3}, it is easy to verify that
\begin{equation}
\left|  \frac{n(\lambda,s)}{d(\lambda,s)}\right| < \rho,
\end{equation}
which corresponds to $\big|G(s)\big|<\rho$. \hfill{\qed}

Following the above developments, here remarkably,
it can be seen that it is possible to note that
$d_{c}(s)$ is appropriate to be interpreted as a central polynomial to characterize the basic performance of the system,
and
additional discussions on the design of an appropriate central polynomial can be found in~\cite{khatibi2008fixed}.

Then, before proceeding to develop the required robust performance criterion, Theorem~\ref{theorem:2} is next readily extended to be suitable for a more general case, which is summarized by the following theorem.

\begin{theorem}~\label{theorem:3}
	Given matrices $Q$, $H_i$, $E_i$, $i=1,2,\cdots,m$ with appropriate dimensions, $Q$ is symmetrical, $\Delta_i=\textup{diag}\{\delta_{i1}, \delta_{i2}, \cdots, \delta_{in}\}$ with $\delta_{ij}\in [-1,1]$, $i=1,2,\cdots, m$, $j=1,2,\cdots,n$,
	holds if and only if there exist matrices $R_i=\textup{diag}\{\varepsilon_{i1}, \varepsilon_{i2}, \cdots, \varepsilon_{in}\}$ with $\varepsilon_{ij}>0$,  $i=1,2,\cdots, m$, $j=1,2,\cdots,n$, such that
	\begin{equation}
	Q+
	\left[
	\begin{array}{ccccc}
	\sum \limits_{i=1}^m E_i^T R_i^{-1}	E_i & 0 & \cdots &0 & 0\\
	0 &  0& \cdots& 0& 0\\
	\vdots& \vdots &  \ddots & \vdots & \vdots\\
	0& 0& \cdots& 0& 0\\
	0&  0& \cdots& 0 &   \sum \limits_{i=1}^m H_i  R_i H_i^T \end{array}
	\right]<0. \label{eqn:theorem3:2}
	\end{equation}
\end{theorem}
\noindent{\textbf{Proof of Theorem~\ref{theorem:3}:}}
Theorem~\ref{theorem:3} is an extension of Theorem~\ref{theorem:2}, and the mostly straightforward proof is omitted. \hfill{\qed}

At this stage, define the transfer functions $G_{p1}(s)$ and $G_{p2}(s)$ as
\begin{equation}\label{DefS1S2}
\begin{aligned}
G_{p1}(s)=&G_{p1n}(s)+G_{p1u}(s),\\
G_{p2}(s)=&G_{p2n}(s)+G_{p2u}(s),
\end{aligned}
\end{equation}
where
\begin{equation}\label{DefS1S2np}
\begin{aligned}
G_{p1n}(s)=&1-G_{sn}(s),\\
G_{p1u}(s)=&-G_{su}(s),\\
G_{p2n}(s)=&\frac{(\bm{a^{c}}\ast\bm{x})\bm{s_{m+n}}^{T}}{d_{c}(s)},\\
G_{p2u}(s)=&\frac{((\bm{a_{d}}\Delta_{a})\ast\bm{x})\bm{s_{m+n-1}}^{T}}{d_{c}(s)}.
\end{aligned}
\end{equation}
From Theorem~\ref{theorem:convex}, $| S(j\omega)|<\rho_s,\, \omega\in\Omega_s$ is guaranteed if the following two conditions hold:\\
\begin{equation}\label{Conv1}
\begin{aligned}
\left|G_{p1}(j\omega)\right|<\delta_s, \quad \omega\in\Omega_s,
\end{aligned}
\end{equation}
and
\begin{equation}\label{Conv2}
\begin{aligned}
\left|G_{p2}(j\omega)\right| <(1-\delta_s)\rho_s, \quad \omega\in\Omega_s,
\end{aligned}
\end{equation}
with $\delta_s\in(0,\,1)$. Since $G_{p1n}(s)$ and $G_{p2n}(s)$ can be realized in the controllable canonical form as
\begin{equation}\label{realization2}
\begin{aligned}
\sum \nolimits_{p1n}&\triangleq
\{A_{p1n}, B_{p1n},	C_{p1n}, D_{p1n}\},\\
\sum \nolimits_{p2n}&\triangleq
\{A_{p2n}, B_{p2n},	C_{p2n}, D_{p2n}\},
\end{aligned}
\end{equation}  respectively, the state-space realizations of $G_{p1}(s)$ and $G_{p2}(s)$ are given by
\begin{equation}
\begin{aligned}
\sum \nolimits_{p1}&\triangleq\{
A_{p1n}, B_{p1n}, C_{p1n}+\bm{a_{d}} \Delta_{a} X+\bm{b_{d}} \Delta_{b} Y, D_{p1n}\},\\
\sum \nolimits_{p2}&\triangleq
\{
A_{p2n}, B_{p2n}, C_{p2n}+\bm{a_{d}} \Delta_{a} X, D_{p2n}
\},
\end{aligned}
\end{equation} respectively.

Then, the conditions of robust performance specification in terms of the sensitivity function are summarized by Theorem~\ref{theorem:performance}.
\begin{theorem}~\label{theorem:performance}
	The robust performance specification $\big| S(j\omega)\big| < \rho_s$, $\omega \in \Omega_s$ of the system~\eqref{PlantModel1} in the presence of bounded parametric uncertainties characterized by standard interval variables is guaranteed under the controller~\eqref{Controller1} if and only if there exist Hermitian matrices $P_{p1}$ and $Q_{p1}>0$, $P_{p2}$ and $Q_{p2}>0$, diagonal matrices $R_{p1a}>0$, $R_{p1b}>0$ and $R_{p2a}>0$, $R_{p2b}>0$ such that
	\begin{equation}\label{LMI2}
	\begin{aligned}
	&\left[
	\begin{array} {c|ccc}
	\Gamma_{p1}
	&
	\begin{array} {c}
	C_{p1n}^T \\
	D_{p1n}^T
	\end{array}
	&\begin{array} {c}
	X^T \\ 0
	\end{array}
	&\begin{array} {c}
	Y^T \\ 0
	\end{array}
	\\ \hline
	\begin{array} {cc}
	C_{p1n} & D_{p1n}
	\end{array}
	&
	H_{p1}
	&0
	&0\\
	\begin{array} {cc}
	X & 0
	\end{array}
	&0
	&-R_{p1a}
	&0\\
	\begin{array} {cc}
	Y & 0
	\end{array}
	&0
	&0
	&-R_{p1b}\\
	\end{array}
	\right]<0,
	\end{aligned}
	\end{equation}
	\begin{equation}\label{LMI3}
	\begin{aligned}
	&\left[
	\begin{array} {c|ccc}
	\Gamma_{p2}
	&
	\begin{array} {c}
	C_{p2n}^T \\
	D_{p2n}^T
	\end{array}
	&\begin{array} {c}
	X^T \\ 0
	\end{array}
	\\ \hline
	\begin{array} {cc}
	C_{p2n} & D_{p2n}
	\end{array}
	&
	H_{p2}
	&0\\
	\begin{array} {cc}
	X & 0
	\end{array}
	&0
	&-R_{p2b}\\
	\end{array}
	\right]<0,
	\end{aligned}
	\end{equation}
	where
	\begin{equation}\label{GKYPDefGamma}
	\begin{aligned}
	\Gamma_{p1}  &=
	\left[
	\begin{array} {cc}
	A_{p1n} & B_{p1n} \\
	I & 0
	\end{array}
	\right]^T
	\Upxi_{p1}
	\left[
	\begin{array} {cc}
	A_{p1n} & B_{p1n}  \\
	I & 0
	\end{array}
	\right]\\&+
	\left[
	\begin{array} {cc}
	0 & 0  \\
	0 & -\delta_s
	\end{array}
	\right],
	\end{aligned}
	\end{equation}
	\begin{equation}
	\begin{aligned}
	\Gamma_{p2} &=
	\left[
	\begin{array} {cc}
	A_{p2n} & B_{p2n} \\
	I & 0
	\end{array}
	\right]^T
	\Upxi_{p2}
	\left[
	\begin{array} {cc}
	A_{p2n} & B_{p2n}  \\
	I & 0
	\end{array}
	\right]\\&+
	\left[
	\begin{array} {cc}
	0 & 0  \\
	0 & (\delta_s-1)\rho_s
	\end{array}
	\right],
	\end{aligned}
	\end{equation}
		$\Upxi_{p1}=\Phi_s\otimes P_{p1} +\Psi_s\otimes Q_{p1}$, $\Upxi_{p2}=\Phi_s\otimes P_{p2}+\Psi_s\otimes Q_{p2}$, $\Phi_s$ and $\Psi_s$ are matrices that characterize the frequency range as given in Section I, $H_{p1}=\bm{a_{d}} R_{p1a} \bm{a_{d}}^T+\bm{b_{d}} R_{p1b} \bm{b_{d}}^T-\delta_s$, $H_{p2}=\bm{a_{d}} R_{p2a} \bm{a_{d}}^T+(\delta_s-1)\rho_s$.
	
\end{theorem}

\noindent{\textbf{Proof of Theorem~\ref{theorem:performance}:}}
From Lemma~\ref{lemma: GKYP}, the robust specification \eqref{Conv1} is guaranteed if and only if there exist Hermitian matrices $P_{p1}$ and $Q_{p1}>0$, such that
\begin{equation}\label{GKYPPerf0}
\begin{aligned}
\Vast[
\begin{array} {c|c}
\Gamma_{p1}
&
\begin{array} {c}
\begin{split} C_{p1n}^T &+X^T \Delta_{a} \bm{a_{d}}^T  \\&+  Y^T \Delta_{b}   \bm{b_{d}}^T \end{split} \\
D_{p1n}^{T}
\end{array}
\\\hline
\begin{array} {cc} \vspace{-3mm}
\begin{split}	C_{p1n}&+\bm{a_{d}} \Delta_{a}X \\&+\bm{b_{d}}  \Delta_{b} Y \end{split}  & D_{p1n}
\end{array}
&
-\delta_s
\end{array}
\Vast]<0
\end{aligned},
\end{equation}


Equivalently, we have
\begin{equation}~\label{GKYPPerf2}
\begin{aligned}
&\left[
\begin{array} {c|c}
\Gamma_{p1}
&
\begin{array} {c}
C_{p1n}^T \\
D_{p1n}^T
\end{array}
\\\hline
\begin{array} {cc}
C_{p1n} & D_{p1n}
\end{array}
&
-\delta_s
\end{array}
\right]\\& +
\left[
\begin{array} {ccc}
0 & 0& X^T \Delta_{a}\bm{a_d}^T  + Y^T \Delta_{b}\bm{b_d}^T\\
0 &	0 & 0\\
\bm{a_{d}}\Delta_{a}X+\bm{b_{d}}\Delta_{b}Y &	0 & 0
\end{array}
\right]
<0.
\end{aligned}
\end{equation}

From Theorem~\ref{theorem:3}, \eqref{GKYPPerf2} holds if and only if there exist positive definite diagonal matrices $R_{p1a}$ and $R_{p1b}$ such that
\begin{equation}\label{GKYPPerf3}
\begin{aligned}
&\left[
\begin{array} {c|c}
\Gamma_{p1}
&
\begin{array} {c}
C_{p1n}^T \\
D_{p1n}^T
\end{array}
\\\hline
\begin{array} {cc}
C_{p1n} & D_{p1n}
\end{array}
&
-\delta_s
\end{array}
\right] \\
&+\left[
\begin{array} {ccc}
X^T R_{p1a}^{-1} X+Y^T R_{p1b}^{-1} Y & 0& 0\\
0 & 0& 0\\
0 & 0& \bm{a_{d}} R_{p1a}\bm{a_{d}}^T+\bm{b_{d}}R_{p1b}\bm{b_{d}}^T
\end{array}
\right] \\&<0.
\end{aligned}
\end{equation}
Then, we have
\begin{equation}\label{GKYPPerf4}
\begin{aligned}
&\left[
\begin{array} {c|c}
\Gamma_{p1}
&
\begin{array} {c}
C_{p1n}^T \\
D_{p1n}^T
\end{array}
\\\hline
\begin{array} {cc}
C_{p1n} & D_{p1n}
\end{array}
&
-\delta_s
\end{array}
\right]
+\left[
\begin{array} {c}
0 \\
0 \\
\bm{a_{d}}
\end{array}
\right]R_{p1a}
\left[
\begin{array} {ccc}
0 & 0 & \bm{a_{d}}^T
\end{array}
\right]\\
&+
\left[
\begin{array} {c}
0 \\
0 \\
\bm{b_{d}}
\end{array}
\right]R_{p1b}
\left[
\begin{array} {ccc}
0 & 0 & \bm{b_{d}}^T
\end{array}
\right]+
\left[
\begin{array} {c}
X^T \\
0 \\
0
\end{array}
\right]R_{p1a}^{-1}
\left[
\begin{array} {ccc}
X & 0 & 0
\end{array}
\right]\\
&+
\left[
\begin{array} {c}
Y^T \\
0 \\
0
\end{array}
\right]R_{p1b}^{-1}
\left[
\begin{array} {ccc}
Y & 0 & 0
\end{array}
\right]<0,
\end{aligned}
\end{equation}
which can be further expressed in the form of~\eqref{LMI2}.

Therefore, it can be obtained that the condition \eqref{Conv1} is guaranteed if and only if there exist Hermitian matrices $P_{p1}$ and $Q_{p1}>0$, and positive definite diagonal matrices $R_{p1a}$ and $R_{p1b}$ such that \eqref{LMI2} is satisfied. Similarly, the condition \eqref{Conv2} is guaranteed if and only if there exist Hermitian matrices $P_{p2}$ and $Q_{p2}>0$, and positive definite diagonal matrices $R_{p2a}$ and $R_{p2b}$ such that \eqref{LMI3} is satisfied. This completes the proof of the theorem. \hfill{\qed}

Furthermore, the results on the specification in terms of the sensitivity function $\big| S(j\omega)\big| < \rho_s$, $\omega \in \Omega_s$ can be extended to the specification on the complementary sensitivity function  $\big| T(j\omega)\big| < \rho_t$, $\omega \in \Omega_t$. Thus define transfer functions $G_{p3}(s)$  as
\begin{equation}
\begin{aligned}
G_{p3}(s)=&G_{p3n}(s)+G_{p3u}(s),
\end{aligned}
\end{equation}
where
\begin{equation}
\begin{aligned}
G_{p3n}(s)=&\frac{(\bm{b^{c}}\ast\bm{y})\bm{s_{m+n}}^{T}}{d_{c}(s)},\\
G_{p3u}(s)=&\frac{((\bm{b_{d}}\Delta_{b})\ast\bm{y})\bm{s_{m+n-1}}^{T}}{d_{c}(s)}.
\end{aligned}
\end{equation}
From Theorem~\ref{theorem:convex}, $| T(j\omega)|<\rho_t,\, \omega\in\Omega_t$ is guaranteed if the following two conditions hold:\\
\begin{equation}
\begin{aligned}
\left|G_{p1}(j\omega)\right|<\delta_t,\quad  \omega\in\Omega_t
\end{aligned},
\end{equation}
and
\begin{equation}
\begin{aligned}
\left|G_{p3}(j\omega)\right| <(1-\delta_t)\rho_t, \quad  \omega\in\Omega_t
\end{aligned},
\end{equation}
with $\delta_t\in(0,1)$. $G_{p3n}(s)$ can be realized in the controllable canonical form as
\begin{equation}\label{realization3}
\begin{aligned}
\sum \nolimits_{p3n}&\triangleq
\{A_{p3n}, B_{p3n},	C_{p3n}, D_{p3n}\},
\end{aligned}
\end{equation} then the state-space realization of $G_{p3}(s)$ is given by
\begin{equation}
\sum \nolimits_{p3} \triangleq
\{
A_{p3n}, B_{p3n}, C_{p3n}+\bm{b_{d}} \Delta_{b} Y, D_{p3n}
\}.
\end{equation}

It is worth mentioning, at this point, that some rather useful properties hold for the state-space realizations of the nominal models \eqref{realization1}, \eqref{realization2}, and \eqref{realization3}, where $A_{sn}=A_{p1n}=A_{p2n}=A_{p3n}$,
$B_{sn}=B_{p1n}=B_{p2n}=B_{p3n}$, $C_{sn}=-C_{p1n}$, $D_{sn}=1-D_{p1n}$.


In a similar manner here, in what follows, Theorem~\ref{theorem:performance2} is proposed to cater to the robust performance specification in terms of the complementary sensitivity function.
\begin{theorem}~\label{theorem:performance2}
	The robust performance specification $\big| T(j\omega)\big| < \rho_t$, $\omega \in \Omega_t$ of the system~\eqref{PlantModel1} in the presence of bounded parametric uncertainties characterized by standard interval variables is guaranteed under the controller~\eqref{Controller1} if and only if there exist Hermitian matrices $P_{p3}$ and $Q_{p3}>0$, $P_{p4}$ and $Q_{p4}>0$, diagonal matrices $R_{p3a}>0$, $R_{p3b}>0$ and $R_{p4a}>0$, $R_{p4b}>0$ such that
	\begin{equation}\label{LMI4}
	\begin{aligned}
	&\left[
	\begin{array} {c|ccc}
	\Gamma_{p3}
	&
	\begin{array} {c}
	C_{p1n}^T \\
	D_{p1n}^T
	\end{array}
	&\begin{array} {c}
	X^T \\ 0
	\end{array}
	&\begin{array} {c}
	Y^T \\ 0
	\end{array}
	\\ \hline
	\begin{array} {cc}
	C_{p1n} & D_{p1n}
	\end{array}
	&
	H_{p3}
	&0
	&0\\
	\begin{array} {cc}
	X & 0
	\end{array}
	&0
	&-R_{p3a}
	&0\\
	\begin{array} {cc}
	Y & 0
	\end{array}
	&0
	&0
	&-R_{p3b}\\
	\end{array}
	\right]<0,
	\end{aligned}
	\end{equation}
	\begin{equation}\label{LMI5}
	\begin{aligned}
	&\left[
	\begin{array} {c|ccc}
	\Gamma_{p4}
	&
	\begin{array} {c}
	C_{p3n}^T \\
	D_{p3n}^T
	\end{array}
	&\begin{array} {c}
	X^T \\ 0
	\end{array}
	\\ \hline
	\begin{array} {cc}
	C_{p3n} & D_{p3n}
	\end{array}
	&
	H_{p4}
	&0\\
	\begin{array} {cc}
	X & 0
	\end{array}
	&0
	&-R_{p4b}\\
	\end{array}
	\right]<0,
	\end{aligned}
	\end{equation}
	where
	\begin{equation}
	\begin{aligned}
	\Gamma_{p3}  &=
	\left[
	\begin{array} {cc}
	A_{p1n} & B_{p1n} \\
	I & 0
	\end{array}
	\right]^T
	\Upxi_{p3}
	\left[
	\begin{array} {cc}
	A_{p1n} & B_{p1n}  \\
	I & 0
	\end{array}
	\right]\\&+
	\left[
	\begin{array} {cc}
	0 & 0  \\
	0 & -\delta_t
	\end{array}
	\right],
	\end{aligned}
	\end{equation}
	\begin{equation}
	\begin{aligned}
	\Gamma_{p4} &=
	\left[
	\begin{array} {cc}
	A_{p3n} & B_{p3n} \\
	I & 0
	\end{array}
	\right]^T
	\Upxi_{p4}
	\left[
	\begin{array} {cc}
	A_{p3n} & B_{p3n}  \\
	I & 0
	\end{array}
	\right]
	\\&+
	\left[
	\begin{array} {cc}
	0 & 0  \\
	0 & (\delta_t-1)\rho_t
	\end{array}
	\right],
	\end{aligned}
	\end{equation}
		$\Upxi_{p3}=\Phi_t\otimes P_{p3} +\Psi_t\otimes Q_{p3}$, $\Upxi_{p4}=\Phi_t\otimes P_{p4}+\Psi_t\otimes Q_{p4}$, $\Phi_t$ and $\Psi_t$ are matrices that characterize the frequency range as given in Section I, $H_{p3}=\bm{a_{d}} R_{p3a} \bm{a_{d}}^T+\bm{b_{d}} R_{p3b} \bm{b_{d}}^T-\delta_t$, $H_{p4}=\bm{b_{d}} R_{p4a} \bm{b_{d}}^T+(\delta_t-1)\rho_t$.
\end{theorem}

\noindent{\textbf{Proof of Theorem~\ref{theorem:performance2}:}}
The proof of Theorem~\ref{theorem:performance2} proceeds along with the same procedures as the proof of Theorem~\ref{theorem:performance}. \hfill{\qed}

Based on all the above developments and analysis,
it can now be noted that with the work in this paper, the
fixed-order controller considering robust stability and performance can be designed by solving
the LMIs \eqref{LMI1}, \eqref{LMI2}, \eqref{LMI3}, \eqref{LMI4}, and \eqref{LMI5}.

\begin{remark}
	For an uncertain system, our proposed methodology formulates the loop shaping problem with fewer conditions, as compared to classical robust control methods. Take an uncertain polytopic system as an example, each vertex is required to be taken into account in the loop shaping problem, and the number of these vertices grows exponentially with the number of parametric uncertainties, which results in significantly increasing conditions in the classical design algorithms. For instance, by using the classical robust control methods, if the number of parametric uncertainties is 8, then a total of 256 vertices are to be considered, which leads to 1280 LMIs to solve.  But with our proposed approach, only one-time checking of matrix existence by solving 5 LMIs is needed exclusively, regardless of the numbers of parameter uncertainties.
\end{remark}

\begin{remark}
	It can also be noted that depending on the design specifications on a specific problem, the conditions on the robust stability, the robust performance in terms of the sensitivity function, and the robust performance in terms of the complementary sensitivity function
	can be implemented either separately or together.
\end{remark}

\begin{remark}
	In closed-loop shaping problem, the specifications on the sensitivity function and the complementary sensitivity function are important indicators to ensure good system performance. The magnitudes of their target values under a specific frequency range are usually pre-defined by the user, which are dependent on the actual problem. It is well known that standard and typical selections of their limiting bounds $\rho_s$ and $\rho_t$ are given by $-3$ dB, and these limiting bounds physically guarantee an effective closed-loop bandwidth and a roll-off frequency for the control system. It is also noted that $\rho_s$ and $\rho_t$ can be simply converted to the optimization objective, so that they can be minimized during the optimization, which physically represents the optimization of disturbance rejection ability and roll off ability, respectively.
\end{remark}

\section{Illustrative example}
Consider the nominal model of a second-order plant with an unstable pole
\begin{equation}
	P(s)=\frac{b_1s+b_2}{s^2+a_1s+a_2},
\end{equation}
with $a_1=3$, $a_2=-10$, $b_1=8$, $b_2=4$. Assume the parametric uncertainties exist in $a_1, a_2$ and  $b_1, b_2$, and the deviations of these parameters are all $\pm$20\% of their nominal values. A second-order controller is to be designed, where
\begin{equation}
	K(s)=\frac{y_0s^2+y_1s+y_2}{s^2+x_1s+x_2}.
\end{equation}
To the industrial preference in terms of practical implementation, $x_2$ is set to be zero so that $K(s)$ becomes a PID controller with a low-pass filter. Here, the typical desirable
design specifications
would
include the robust stability and the robust performance with the sensitivity function $\big| S(j\omega)\big|<-3$ {dB} under the frequency range $\omega \in [$0.01 rad/s, 0.1 rad/s$]$, and the complementary sensitivity function $\big| T(j\omega)\big|<-3$ {dB} under the frequency range $\omega \in [$20 rad/s, 100 rad/s$]$ (essentially stated in this equivalent manner).

\begin{figure}[t]
	\centering
	\includegraphics[trim=0 0 0 0,width=1\columnwidth]{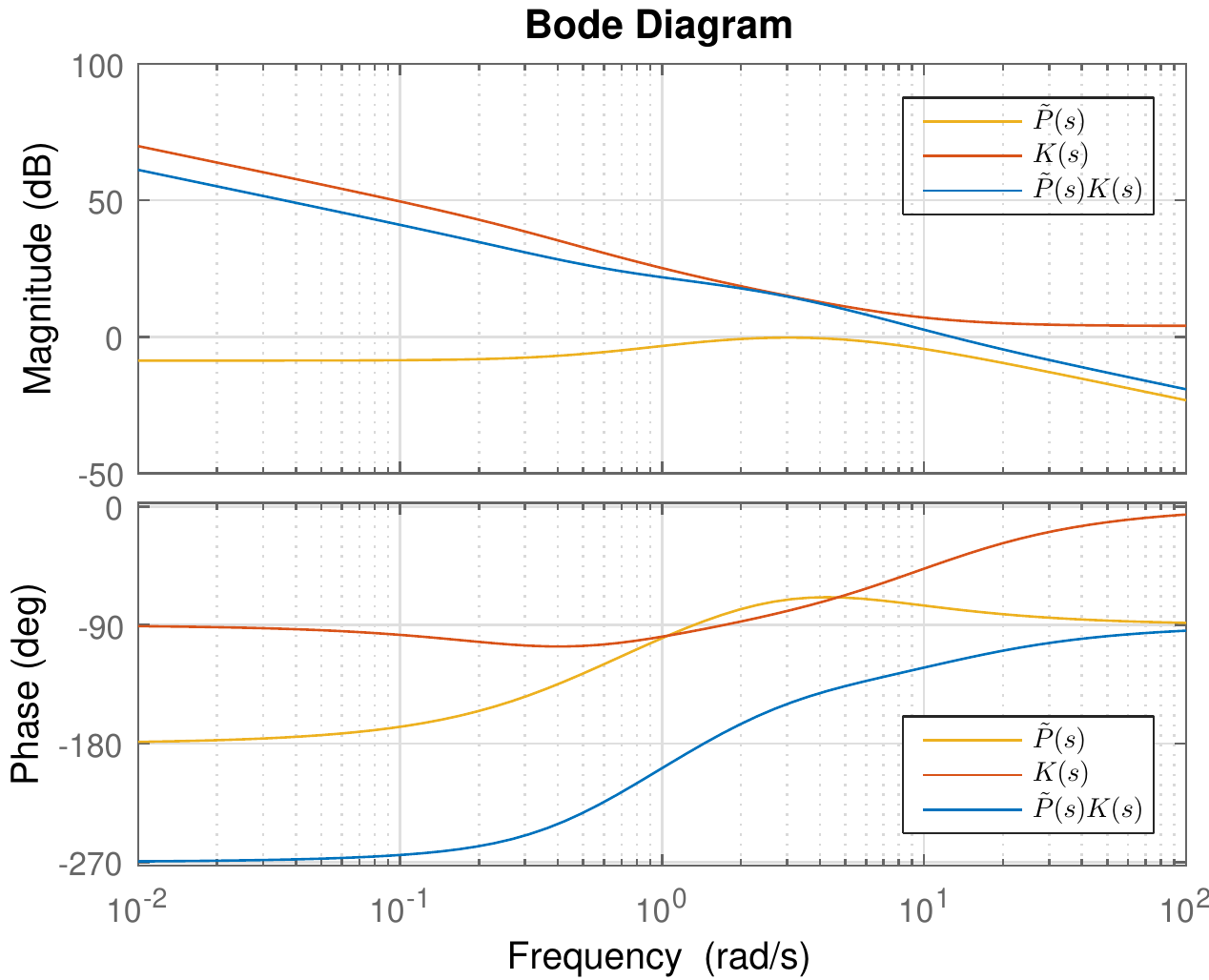}
	\caption{Bode diagrams of the plant, the controller, and the open-loop system.}
	\label{fig:openloop}
\end{figure}

\begin{figure}[t]
	\centering
	\includegraphics[trim=0 0 0 0,width=1\columnwidth]{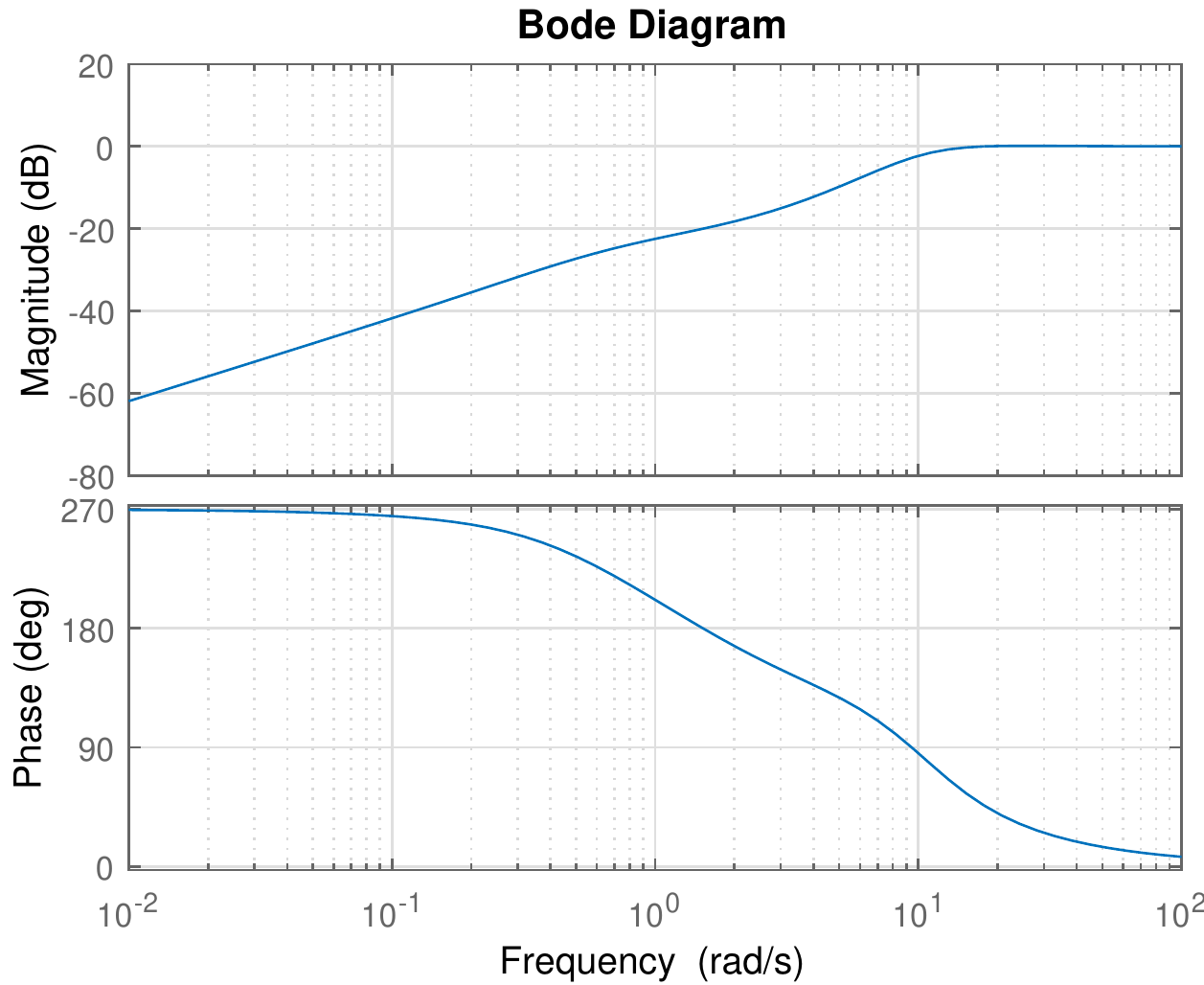}
	\caption{Bode diagram of the sensitivity function.}
	\label{fig:sensitivity}
\end{figure}

\begin{figure}[t]
	\centering
	\includegraphics[trim=0 0 0 0,width=1\columnwidth]{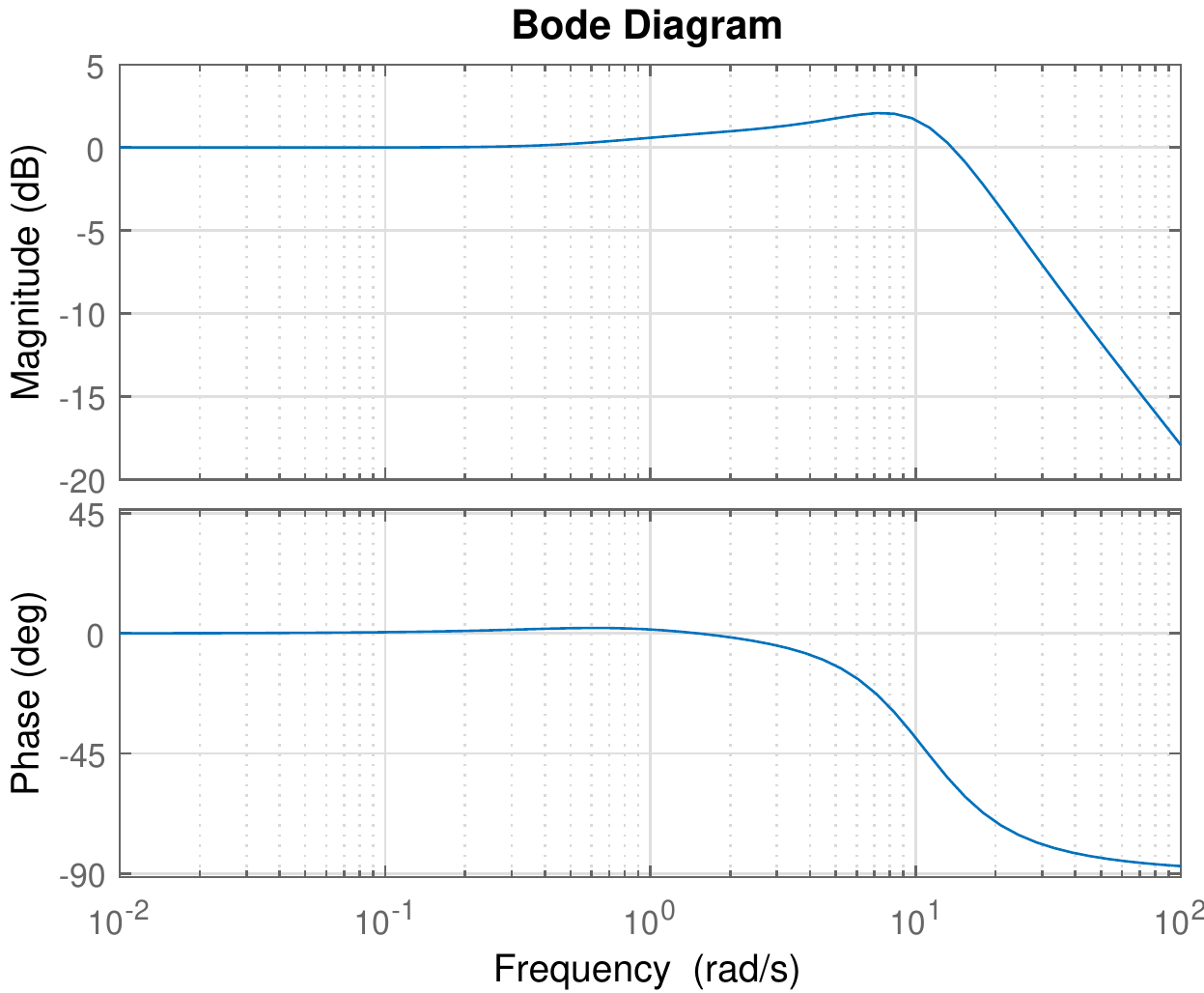}
	\caption{Bode diagram of the complementary sensitivity function.}
	\label{fig:csensitivity}
\end{figure}

With the methodology proposed here, thus
the design can proceed with the steps where
a Hurwitz polynomial is chosen as $d_c(s)= s^4+16s^3+89s^2+390s+200$ and the parameters $\delta_s$ and $\delta_t$ are chosen as 0.5.
The LMIs \eqref{LMI1}, \eqref{LMI2}, \eqref{LMI3}, \eqref{LMI4}, and \eqref{LMI5} are
then constructed and solved by the YALMIP Toolbox in MATLAB.
The resulting controller parameters are given by $x_1=0.3307$, $y_0=1.5790$, $y_1=16.9886$, $y_2=10.2572$.
In the simulation, we use a plant $\tilde P(s)$ with one set of the uncertain parameters given by $\tilde a_1=3.5486$, $\tilde a_2=-9.9415$, $\tilde b_1= 6.9044$, $\tilde b_2=3.6471$. The controller parameters obtained by the proposed approach are then implemented on the uncertain plant $\tilde P(s)$, and it can be easily verified that the closed-loop system is stabilized, where the closed-loop poles are $-6.7471 \pm 7.1417i$, $-$0.8072, and $-$0.4801. The Bode diagrams of the uncertain plant $\tilde P(s)$, the controller $K(s)$, and the open-loop system $\tilde P(s)K(s)$ are shown in Figure~\ref{fig:openloop}. Also, the Bode diagrams of the sensitivity function $S(j\omega)$ and the complementary sensitivity function $T(j\omega)$ are illustrated in Figure~\ref{fig:sensitivity} and Figure~\ref{fig:csensitivity}, respectively. It is shown that all the design specifications are met by using the proposed approach, and the stated
robust stability and performance are successfully achieved with a fixed-order controller design
under this situation of an uncertain system.


In the following work, the case by using our proposed method is denoted by Case I. For comparative purposes, the fixed-order controller design approach adopted in~\cite{zhu2019dual} is used with the same Hurwitz polynomial as our proposed approach. With this approach, first of all, all the uncertain systems are considered, and this case is denoted by Case II. Since there are 4 uncertain parameters in the plant, the number of the resulting uncertain systems is given by $2^4=16$. As a result, 80 LMIs are integrated into the program considering robust stability and robust performance in terms of the sensitivity function and the complementary sensitivity function. In this case, the controller parameters are given by $x_1= 7.2850$, $y_0=1.0875$, $y_1=13.1769$, $y_2=68.6918$.  Secondly, only the nominal system is considered in this method, and this case is denoted by Case III. In this case, only 5 LMIs are required in the program, and the controller parameters are given by $x_1=12.7610$, $y_0=0.4414$, $y_1=11.6312$, $y_2=62.9964$.

In terms of the computational effort, the computational time in all three cases is recorded and summarized in Table~\ref{table:tab}. As can be seen, the computational time in Case I and Case III is less than half of Case II, and the main reason is that Case II considers the whole uncertain domain by checking all the 4 uncertain systems. In fact, the number of uncertain systems grows exponentially with the number of parametric uncertainties in the plant. When the number of uncertainties increases, the computation is rather costly and the difference in computational time between Case I and Case II would be more significant and more clearly observed. Compared with Case II, our proposed methodology manipulates all these parametric uncertainties as a whole, and thus only a one-time checking of matrix existence by solving 5 LMIs is needed exclusively, regardless of the number of parameter uncertainties. Moreover, the difference in computational time between Case I and Case III is minimal. However, Case III only takes the nominal system into consideration, and thus the system performance would be degraded if the system is perturbed.

\begin{figure}[t]
	\centering
	\includegraphics[trim=0 0 0 0,width=1\columnwidth]{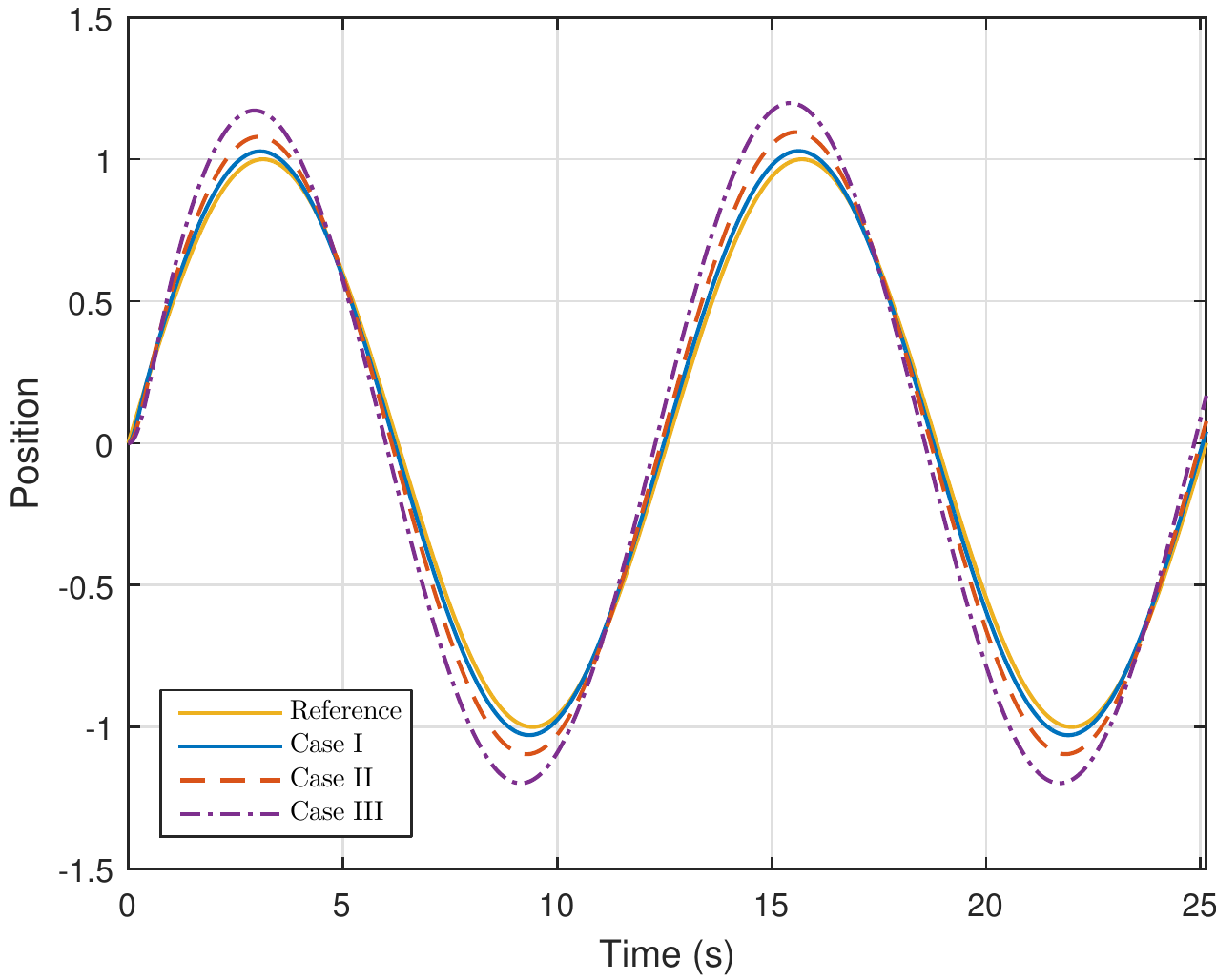}
	\caption{Comparison of tracking performance.}
	\label{fig:sinetrack}
\end{figure}

To clearly observe the comparison of these methods in terms of system performance, simulations in the time domain are conducted, and the objective is to track a sinusoidal signal. In the simulation, the uncertain plant $\tilde P(s)$ is used again, and Figure~\ref{fig:sinetrack} depicts the comparison of tracking performance in the closed-loop control with different methods. Besides, the resulting root-mean-square error (RMSE) and the maximum absolute error (MaxAE) is given in Table~\ref{table:tab}. It is straightforward to see that Case I achieves the best system performance in this task. On the other hand, the performance in Case III is the worst. The main reason behind this phenomenon is that Case III does not take the parametric uncertainties into account, and thus the controller parameters do not have enough robustness to resist the effects of parametric uncertainties. Since the plant used in this simulation is perturbed, the system performance is downgraded.

To sum up, the proposed approach has demonstrated its effectiveness and efficiency in solving a class of loop shaping problems. Rather importantly, the proposed method maintains good robustness towards parametric uncertainties, and also avoids the computational burden resulted from the necessity of checking the robust stability and performance for all uncertain systems in the traditional robust control techniques.

\begin{table}[t!]\centering
	\caption{Comparison of computational efficiency and system performance}
	\label{table:tab}
	\begin{tabular}{|c|c|c|c|}
		\hline
		& Case I & Case II & Case III \\ \hline
		Time (s) & 1.0883 & 2.6000  & 1.1206   \\ \hline
		RMSE     & 0.0325 & 0.0794  & 0.1652   \\ \hline
		MaxAE    & 0.0488 & 0.1205  & 0.2498   \\ \hline
	\end{tabular}
\end{table}

\section{Conclusion}
In this work, a fixed-order robust controller design approach is developed under a specific restricted frequency range.
To achieve this, first of all,
an initial set of newly developed theoretical results to be used in the robust stability and robust performance criteria are presented.
Secondly, the robust stabilization condition for the
situation of uncertain systems is constructed with the concept of positive realness.
Thirdly, the robust performance specifications are characterized under a restricted frequency range;
and the frequency-domain system performance from the viewpoint of sensitivity shaping
is realized in a time-domain framework.
These conditions for robust stability and robust performance are given
and formulated by the respective LMIs.
An illustrative example of a relevant appropriate controller design problem is given
and the effectiveness and efficiency of the proposed theoretical results are validated from the comparative simulations.

\bibliographystyle{IEEEtran}
\bibliography{IEEEabrv,Reference}

\end{document}